\documentclass[12pt]{article}
\usepackage[final]{epsfig}
\usepackage{graphics}
\usepackage{amsmath}
\usepackage{amsfonts}
\usepackage{latexsym}
\usepackage{amssymb}
\usepackage{graphicx}
\usepackage{epstopdf}
\newtheorem{lemma}{Lemma}[section]

\newtheorem{remark}[lemma]{Remark}

\newtheorem{theorem}{Theorem}

\begin{document}
\newcommand{\eps}{{\varepsilon}}
\newcommand{\proofend}{$\Box$\bigskip}
\newcommand{\C}{{\mathbf C}}
\newcommand{\Q}{{\mathbf Q}}
\newcommand{\R}{{\mathbf R}}
\newcommand{\Z}{{\mathbf Z}}
\newcommand{\RP}{{\mathbf {RP}}}

\newcommand\pd[2]{\frac{\partial #1}{\partial #2}}
\def\proof{\paragraph{Proof.}}

\title {Variations on the Tait--Kneser  theorem}
\author{Serge Tabachnikov and Vladlen Timorin\\
{\it Department of Mathematics, Penn State University}\\
{\it University Park, PA 16802, USA}\\
{\it Institute for Mathematical Sciences,}\\
{\it State University of New York at Stony Brook}\\
{\it Stony Brook, NY  11794, USA}
}
\date{}
\maketitle

\section{Introduction} \label{intro}

At every point, a smooth plane curve can be approximated, to
second order, by a circle; this circle is called osculating. One
may think of the osculating circle as passing through three
infinitesimally close points of the curve. A {\em vertex} of the
curve is a point at which the osculating circle hyper-osculates:
it approximates the curve to third order. Equivalently, a vertex
is a critical point of the curvature function.

\begin{figure}[hbtp]
\centering
\includegraphics[width=3in]{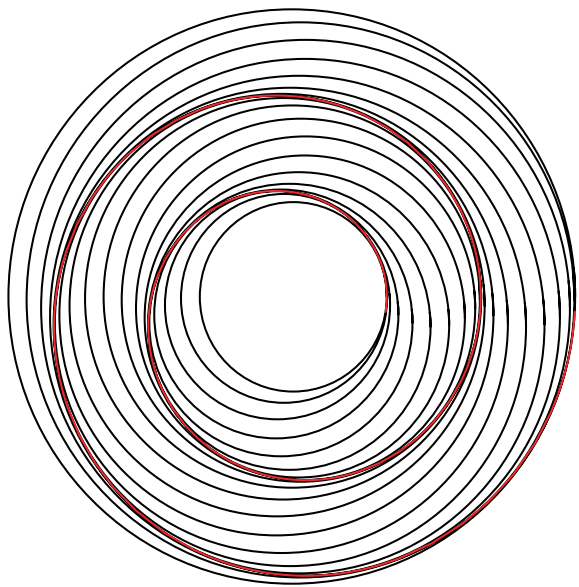}
\caption{A spiral and its nested osculating circles}
\label{circ}
\end{figure}

Consider a (necessarily non-closed) curve, free from vertices. The
classical Tait-Kneser theorem \cite{Tai,Kne}  (see also \cite{Gug,O-T}),
states that the osculating circles of the curve are pairwise
disjoint, see Figure \ref{circ}. This theorem is closely related
to the four vertex theorem of S. Mukhopadhyaya \cite {Muh}  that a
plane oval has at least 4 vertices (see again \cite{Gug,O-T}).

Figure \ref{circ} illustrates the Tait-Kneser theorem: it shows an
annulus foliated by osculating circles of a curve.

\begin{remark} \label{nondif}
{\rm This  foliation is not differentiable! Here is a proof. Let
$f$ be a differentiable function in the annulus, constant on the
leaves. We claim that $f$ is constant. Indeed, $df$  vanishes on
the tangent vectors to the leaves. The curve  is tangent to its
osculating circle at every point, hence $df$ vanishes on the curve
as well. Hence $f$ is constant on the curve. But the curve
intersects all the circles that form the annulus, so  $f$ is
constant everywhere.}
\end{remark}

\begin{remark}
{\rm The Tait-Kneser theorem has an  analog in plane Minkowski geometry,
see \cite{Tab}.}
\end{remark}

We will prove a number of analogs of the Tait-Kneser theorem; in each
case, we will obtain a non-differentiable foliation with smooth
leaves.

\section{Osculating Taylor polynomials} \label{Tpoly}

Let $f$ be a smooth function of one real variable. Fix $n \geq 1$
and let $t\in \R$. The osculating (Taylor) polynomial $g_t$ of
degree $n$ of the function $f$ at the point $t$ is the polynomial,
whose value and the values of whose first $n$ derivatives at the
point $t$ coincide with those of $f$:
\begin{equation}
\label{Taylor}
g_t (x) =\sum_{i=0}^n \frac{f^{(i)}(t)}{i!} (x-t)^i.
\end{equation}
The osculating polynomial $g_t$ is {\em hyper-osculating} if it
approximates the function $f$ at the point $t$ up to $n+1$-st
derivative, that is, if $f^{(n+1)}(t)=0$.

Assume that $n$ is even and $f^{(n+1)} (t)\neq0$ on some interval
$I$ (possibly, infinite).

\begin{theorem} \label{Knpol}
For any distinct $a,b \in I$, the graphs of the osculating
polynomials $g_a$ and $g_b$ are disjoint.\footnote{We will not be surprised if this result is not new but we did not see it in the literature.}

\end{theorem}

\paragraph{Proof.}
To fix ideas, assume that $f^{(n+1)} (t)>0$ on $I$. Let $a<b$ and
suppose that $g_a(x)=g_b(x)$ for some $x\in\R$. It follows
from (\ref{Taylor}) that
$$
\frac{\partial g_t}{\partial t} (x)= \sum_{i=0}^n \frac{f^{(i+1)}(t)}{i!} (x-t)^i - \sum_{i=0}^n \frac{f^{(i)}(t)}{(i-1)!} (x-t)^{i-1} = \frac{f^{(n+1)}(t)}{n!} (x-t)^n,
$$
and hence $(\partial g_t/\partial t) (x)> 0$ (except for $t=x$). It follows that $g_t(x)$ increases, as a function of $t$, therefore  $g_a(x)<g_b(x)$. This is a contradiction.
\proofend

The same argument proves the following variant of Theorem
\ref{Knpol}. Let $n$ be odd. Assume that  $f^{(n+1)} (t)\neq0$ on
an interval $I$. Consider two points  $a<b$ from $I$.

\begin{theorem} \label{Knpol1}
The graphs of the osculating polynomials $g_a$ and $g_b$ are
disjoint over the segment $[b,\infty)$.
\end{theorem}

\begin{figure}[hbtp]
\centering
\includegraphics[width=4in]{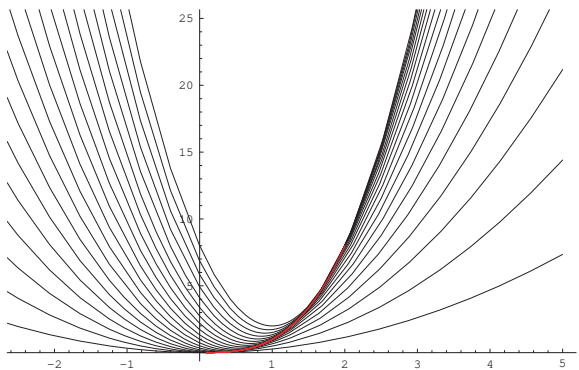}
\caption{Osculating quadratic polynomials of the function $f(x)=x^3$}
\label{quadr}
\end{figure}

\begin{figure}[hbtp]
\centering
\includegraphics[width=4in]{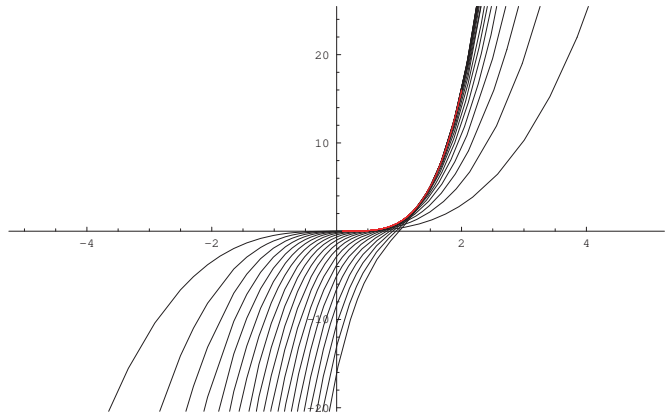}
\caption{Osculating cubic polynomials of the function $f(x)=x^4$}
\label{cube}
\end{figure}

Figure \ref{quadr} shows the graphs of the osculating quadratic polynomials of the function
$f(x)=x^3$ and Figure \ref{cube} of the osculating cubic polynomials of the function $f(x)=x^4$.

\section{Osculating trigonometric polynomials} \label{trpoly}

Let $f$ be a $2\pi$-periodic smooth function, that is, a function
on the circle $S^1=\R/2\pi \Z$. Fix $n \geq 1$ and let $t\in S^1$.
A trigonometric polynomial of degree $n$
$$
g_t (x)=c+\sum_{i=1}^n (a_i \cos ix + b_i \sin ix)
$$
is the {\em osculating trigonometric polynomial} of the function
$f(x)$ at the point $t$ if its value and the values of its first
$2n$ derivatives at the point $t$ coincide with those of $f$.

\begin{remark}
{\rm
The osculating trigonometric polynomial always exists.
Actually, the following more general fact is classical (see, e.g., \cite{Th-U}).
Let $f_i$, $i=1,\dots,N$, be a system of functions on an interval $I$ such that
the Wronski determinant of this system is nonzero everywhere on $I$.
Then, for any sufficiently smooth function $g$ on $I$ and any $t_0\in I$,
there is a linear combination of functions $f_i$ that, at $t_0$, approximates $g$ up to
the derivative of order $N-1$.
This boils down to solving the linear system
$$
g^{(j)}(t_0)=\sum c_i f_i^{(j)}(t_0),\quad j=0,\dots,N-1
$$
with unknowns $c_i$, which has a solution due to non-zero determinant.
The solution depends smoothly on $t_0$.

In our case, the functions $f_i$ are $1$, $\cos t$, $\sin t$, $\cos 2t$,
$\sin 2t$, etc., and $N=2n+1$.
The Wronskian of these functions is constant, which can be seen by
differentiating its columns.
On the other hand, the functions are linearly independent solutions
of a $N$-th order linear differential equation, hence
the Wronskian is nonzero.

Geometrically, we consider $N$-dimensional projective space and the curve $[f_1:...:f_N:g]$.
The osculating hyperplane of this curve at the point
$t_0$ approximates the curve with $N-1$ derivatives.
The equation of this hyperplane is $g= \sum c_i f_i$,
and this gives the desired approximation.
}
\end{remark}

The osculating trigonometric polynomial $g_t$ is {\em
hyper-osculating} if it approximates the function $f$ at the point
$t$ up to $2n+1$-st derivative, that is, if
$f^{(2n+1)}(t)=g_t^{(2n+1)}(t)$. Trigonometric polynomials of
degree $n$ are annihilated by the differential operator ${\cal
D}:=d(d^2+1)(d^2+4)\dots (d^2+n^2)$, where $d=d/dx$. Therefore
$g_t$ hyper-osculates a function $f$ if and only if $({\cal D}f)
(t)=0$.

Assume that the osculating trigonometric polynomials of degree $n$
for a function $f$ do not hyper-osculate on an interval $I \subset
S^1$.

\begin{theorem} \label{trigThm}
For any distinct $a,b \in I$, the graphs of the osculating
trigonometric polynomials $g_a$ and $g_b$ are disjoint.
\end{theorem}

\paragraph{Proof.}
It is not hard to see that the real number $g_t^{(2n+1)}(t)$
depends continuously on $t$ (indeed, the function $g_t$ depends
continuously on $t$ in the $C^{2n+1}$-metric).

To fix ideas, assume that $f^{(2n+1)}(t)>g_t^{(2n+1)} (t)$ for all $t\in I$.
We will show that $\partial g_t(x)/\partial t > 0$ for all $t\in I$ and all $x\in
S^1$ (except $t=x$), and this will imply the statement of the
theorem as in the proof of Theorem \ref{Knpol}.

Since $g_t$ is an osculating trigonometric polynomial, one has:
\begin{equation} \label{osceq}
g_t^{(j)} (t) = f^{(j)} (t),\quad j=0,\dots,2n.
\end{equation}
Differentiate:
$$
\frac{\partial g_t}{\partial t}^{(j)} (t) + g_t^{(j+1)} (t) = f^{(j+1)} (t),
$$
and combine with (\ref{osceq}) to obtain:
\begin{equation} \label{zeros}
\frac{\partial g_t}{\partial t}^{(j)} (t)=0,\ \ j=0,\dots,2n-1;\quad \frac{\partial g_t}{\partial t}^{(2n)} (t) + g_t^{(2n+1)} (t) = f^{(2n+1)} (t).
\end{equation}

The function $\partial g_t/\partial t$ is a trigonometric polynomial of degree $n$.
If this trigonometric polynomial is not identically zero,
then it has no more than $2n$ roots, counting with multiplicities.
If $\partial g_t/\partial t\equiv 0$, then $(\partial g_t/\partial t)^{(2n)} (t) =0$,
and the last equality in (\ref{zeros}) implies that $g_t$ hyper-osculates.
Thus $\partial g_t/\partial t$ is not identically zero.

 According to (\ref{zeros}), the trigonometric polynomial
$\partial g_t/\partial t$ already has a root at the point $t$ of multiplicity $2n$.
Hence $(\partial g_t/\partial t)(x) \neq 0$ for $x \neq t$.
By the assumption made at the beginning of the proof and the last equality in (\ref{zeros}),
we have $(\partial g_t/\partial t)^{(2n)} (t) >0$.
Hence $(\partial g_t/\partial t)(x) >0$ for $x$ sufficiently close to $t$, and therefore
$(\partial g_t/\partial t)(x) >0$ for all $x\neq t$.
\proofend

\begin{figure}[hbtp]
\centering
\includegraphics[width=4in]{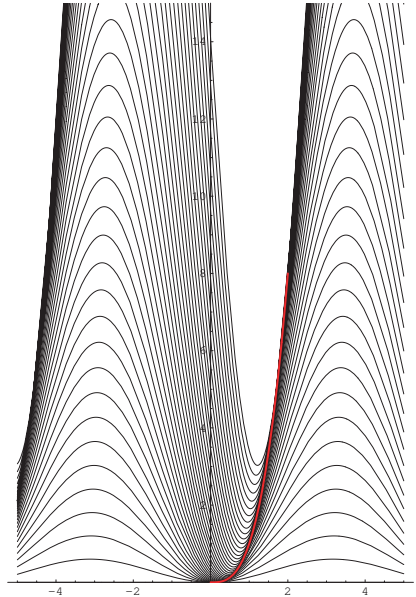}
\caption{Osculating linear harmonics of the function $f(x)=x^3$}
\label{trig}
\end{figure}

Theorem \ref{trigThm} is illustrated in Figure \ref{trig} depicting the graphs of
osculating linear harmonics $c+a\cos x + b \sin x$ for the function $f(x)=x^3$.

\begin{remark} \label{Cheb}
{\rm Theorem \ref{trigThm} extends from trigonometric polynomials
to Chebyshev systems of functions; the proof remains the same.}
\end{remark}

\begin{remark}
{\rm For $n=1$, Theorem \ref{trigThm} implies the Tait-Kneser theorem: it suffices to
consider the support function of the curve and use the fact that the support
functions of circles are linear harmonics.}
\end{remark}

\section{Osculating conics, cubics and fractional
linear transformations} \label{cubic}

Fix $d\geq 1$ and consider the space of algebraic curves of degree $d$.
This space has dimension $n(d)=d(d+3)/2$.
At every point, a smooth plane curve $\gamma$ can be approximated, to order $n(d)-1$,
by an algebraic curve of degree $d$;
this algebraic curve is called the {\em osculating curve}.
One may think of the osculating algebraic curve as passing through $n(d)$
infinitesimally close points of $\gamma$.
A {\em $d$-extactic point} of the curve $\gamma$ is a point, at which the osculating
algebraic curve {\em hyper-osculates}: it approximates $\gamma$ to order $n(d)$; see \cite{Ar}.

In this section, we extend the Tait-Kneser theorem to osculating conics and osculating cubic curves.
We assume that the curve $\gamma$ is free from extactic points.
We also assume that the osculating conics and cubic curves along $\gamma$ are non-degenerate.

Consider a smooth function $f$ with nowhere vanishing derivative.
For every  $t\in \R$, there exists a fractional-linear transformation $g_t$,
whose value and the value of whose first two derivatives at the point $t$ coincide
with those of $f$; this is the {\em osculating fractional-linear transformation}.
As before, it {\em hyper-osculates} at the point $t$ if the third derivatives coincide as well.
This happens if and only if the Schwarzian derivative of $f$ vanishes:
$$
\left( \frac{f'''}{f'}- \frac{3}{2} \left( \frac{f''}{f'} \right)^2 \right) (t)=0.
$$
The graph of a fractional-linear transformation is a hyperbola with vertical and
horizontal asymptotes (or a straight line); we refer to these graphs as the
{\em osculating hyperbolas}.
Assume that the osculating hyperbolas for the function $f$
do not hyper-osculate on an interval $I$.
Let $\gamma$ be the graph of $f$ over $I$.

From the projective point of view, all non-degenerate conics are equivalent;
since our results are projectively-invariant,
we assume that the osculating conics of $\gamma$ are ellipses.
In the case of cubic curves, we assume that the osculating cubics of $\gamma$ have two components,
an oval and a branch going to infinity, and that  the ovals, not the infinite branches, osculate $\gamma$.
With these assumptions, we have the next theorem.

\begin{theorem} \label{twothree}
1) The osculating ellipses along $\gamma$ are pairwise disjoint;\\
2) the ovals of the osculating cubic curves along $\gamma$ are pairwise disjoint;\\
3) the osculating hyperbolas of $\gamma$ are pairwise disjoint.
\end{theorem}

\begin{figure}[hbtp]
\centering
\includegraphics[width=5in]{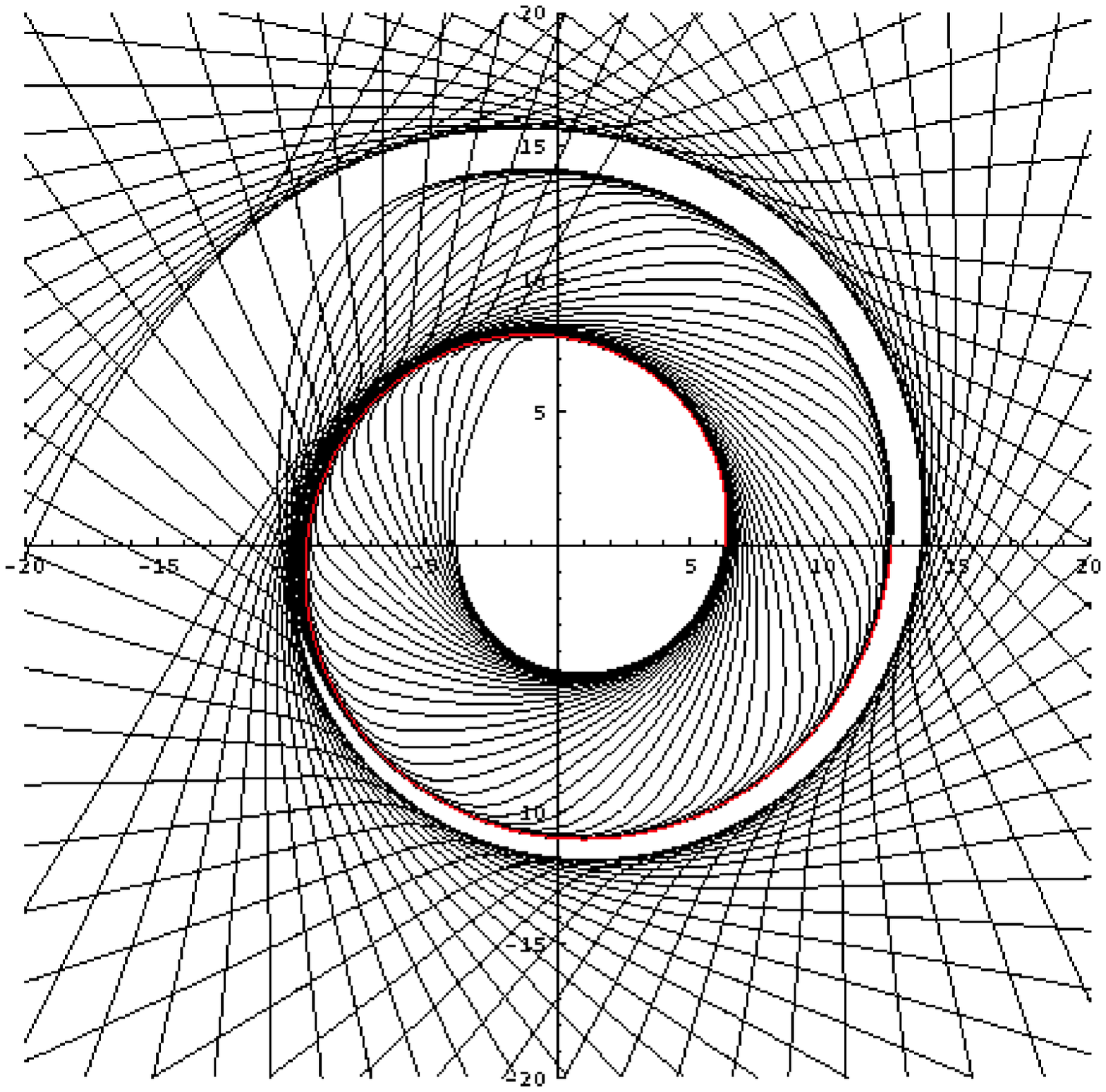}
\caption{Osculating cubic curves of a spiral}
\label{cubics}
\end{figure}

Theorem \ref{twothree}, case 2) is illustrated by Figure \ref{cubics}.

\paragraph{Proof of Theorem \ref{twothree}.}
We will argue about cubic curves, indicating the difference with the case of conics
and hyperbolas, when necessary.

Give the curve $\gamma$ a smooth parameterization, $\gamma(s)$.
Let $\Gamma_s$ be the oval of the osculating cubic curve at the point $\gamma(s)$,
and $f_s (x,y)=0$ its cubic equation.
It suffices to prove that the  curves $\Gamma_a$ and $\Gamma_b$ are nested for distinct
parameter values $a$ and $b$, sufficiently close to each other.

Give the ovals $\Gamma_s$ a smooth parameterization, $\Gamma_s(t)$, such that
the tangency point with the curve $\gamma$ corresponds to $t=0$, that is, $\Gamma_s(0)=\gamma(s)$.
Let $F$ be the map $(s,t) \mapsto \Gamma_s(t)$.
We claim that, for $t \neq 0$, this map is an immersion.
This claim implies that $\Gamma_a$ and $\Gamma_b$ are nested for sufficiently close $a$ and $b$.

To prove the claim, we need the following lemma.

\begin{lemma} \label{Jac}
Suppose that $F(s,t)=(x,y)$. The Jacobian of $F$ vanishes at point $(s,t)$ if and only if
\begin{equation} \label{syst}
\frac{\partial f_s}{\partial s}(x,y)=0,\quad f_s(x,y)=0.
\end{equation}
\end{lemma}

\paragraph{Proof of Lemma.}
The covector $df_s$ is nowhere zero since the curve $\Gamma_s$ is non-degenerate.
This covector vanishes on $\partial F/\partial t$, the tangent vector to the curve $\Gamma_s$.
Therefore the Jacobian of $F$ vanishes exactly when $df_s$ also vanishes on the vector $\partial F/\partial s$.

Differentiate the equation $f_s(F)=0$ with respect to $s$:
$$
\frac{\partial f_s}{\partial s}\circ F+df_s\left(\frac{\partial F}{\partial s}\right)=0.
$$
Thus $df_s$ vanishes on the vector $\partial F/\partial s$ if and only if
$f_s=0$ and $\partial{f_s}/\partial s=0$.
\proofend

Now we need to prove that the system of equations (\ref{syst}) has no
solutions for $t\neq 0$ and point $(x,y)$ on the oval $\Gamma_s$.
Both equations in (\ref{syst}) are cubic, and they are not proportional
since $\gamma(s)$ is not an extactic point (in the case of osculating conics,
the two equations are quadratic).
By the Bezout theorem, the number of solutions is at most 9 (and 4, for conics).
In the case of hyperbolas, $f_s (x,y)=(x-a)(y-b)-c$ where $a,b$ and $c$ depend on $s$;
hence $\partial{f_s}/\partial s=0$ is a linear equation in $x$ and $y$, and system (\ref{syst})
has at most 2 solutions.

For any parameter value $s$, the point $\gamma(s)$ is a multiple solution of system (\ref{syst}).
Since the curve $\{f_s=0\}$ is the osculating curve of degree $d$ for the
curve $\gamma$ at the point $\gamma(s)$, the function $s'\mapsto f_s(\gamma(s'))$
has zero of order $n(d)$ at point $s'=s$.
We can view $f_s(\gamma(s'))$ as a smooth function of two variables $s$ and $s'$.
This function vanishes on the line $s'=s$.
According to a version of the preparation theorem for differentiable functions
\cite{Hor,Loj} (see also \cite{Mal}),
there exists a smooth function $\phi$ of two variables such that
$$
f_{s}(\gamma(s'))=(s-s')^m\phi(s,s')
$$
and $\phi(s,s)\ne 0$ locally near a given value of $s$.
Restricting this equation to a line $s=const$, we obtain $m=n(d)$.
Differentiating with respect to $s$, we see
that $\pd{f_s} s(\gamma(s'))$ starts with terms of order $n(d)-1$ in $s-s'$.
Then $\pd{f_s}s(\Gamma_s(t))$ vanishes for $t=0$ with order $n(d)-1$,
because $\Gamma_s$ approximates $\gamma$ up to order $n(d)$ at $\gamma(s)$.
Hence the multiplicity of the solution $\gamma(s)$ of system (\ref{syst})
is $n(d)-1$.

For $d=2$ (the case of osculating ellipses), this multiplicity is 4, and hence there are no other solutions.
For $d=3$ (the case of osculating cubics), the multiplicity is 8, and there may be one other solution.
However, the number of intersection points of an oval with any  curve is even,
and therefore the 9-th point (if it exists) lies on the other branch.
Therefore system (\ref{syst}) has no  solutions for $t\neq 0$.

Finally, in the case of hyperbolas, the multiplicity of the solution of system (\ref{syst}) at
the point $\gamma(s)$ is 2, therefore there are no other solutions again.
This completes the proof.
\proofend

\begin{remark} \label{Arn}
{\rm It is interesting to compare Theorem \ref{twothree} with three
results on the existence of ``vertices": a  plane oval has at least six sextactic (i.e., 2-extactic)
points \cite{Muh}; a closed plane curve, sufficiently close to an oval of a cubic curve,
has at least ten 3-extactic points \cite{Ar}; and the Schwarzian derivative of a diffeomorphism of $\RP^1$ has at least four zeros \cite{Ghy} (see also \cite{O-T1,O-T}). Likewise, Theorem \ref{trigThm} is related to the existence of ``flexes", i.e., hyperosculating trigonometric polynomials \cite{Th-U}.
}
\end{remark}

\begin{remark}
{\rm In fact, the osculating hyperbolas {\it are} the osculating circles in
Lorentz metric \cite{Ghy,Tab1}.}
\end{remark}

\begin{figure}[hbtp]
\centering
\includegraphics[width=3in]{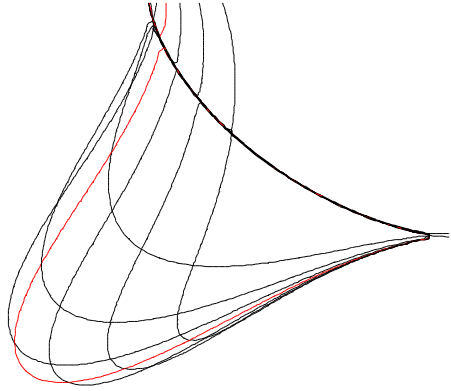}
\caption{Osculating quartics of the curve $x^{2/3}+y^{2/3}=1$.}
\label{osc_quartics}
\end{figure}

\begin{remark}
{\rm
Theorem \ref{twothree} does not generalize to osculating quartics.
This can be seen on Figure \ref{osc_quartics}, where several osculating quartics
for the curve $x^{2/3}+y^{2/3}=1$ are drawn.
Each quartic in the picture splits into two ovals, one being
below and one above the curve.
One can see that nearby ovals below the curve intersect.
}
\end{remark}

\section{Infinitesimal intersection indices}

In this section, we give some more general results that may highlight the proof
of Theorem \ref{twothree}.

Consider a smooth map $F$ of a region in $\R^2$ to a region in $\R^2$.
The map $F$ gives rise to a family of curves.
Namely, for any $s\in\R$, we have the parameterized curve $\Gamma_s:t\mapsto F(s,t)$,
where the parameter $t$ runs through all real numbers such that $(s,t)$ is in the
domain of $F$.
Suppose that the curve $\Gamma_s$ is given locally by an equation $f_s=0$,
which depends smoothly on $s$.
We will assume that $df_s$ never vanishes (e.g., if $f_s$ are polynomials,
then we are talking about nonsingular algebraic curves $\Gamma_s$).

Let $(x,y)$ be a point $\Gamma_s(t)$ on a curve $\Gamma_s$ so that $F(s,t)=(x,y)$.
Define the {\em infinitesimal intersection multiplicity} of $\Gamma_s$ at point
$(x,y)$ as the order of vanishing of the function
$$
t\mapsto Jacobian[F](s,t)
$$
at point $t$.
In particular, if the infinitesimal intersection multiplicity is zero,
then the family $F$ looks like a foliation locally near the point $(x,y)$
and for parameter values near $s$ (however, the curves from the family $F$
corresponding to far-away parameter values may also pass through $(x,y)$).
The {\em infinitesimal intersection index} of a curve $\Gamma_s$ (in the family $F$)
is the sum of local intersection multiplicities at all points of this curve.
The following theorem is an infinitesimal version of the classical Bezout
theorem:

\begin{theorem}
\label{Bez}
Suppose that all curves $\Gamma_s$ are algebraic of degree $d$.
Then the infinitesimal intersection index of each curve $\Gamma_s$ is at most $d^2$.
\end{theorem}

The proof of this theorem is based on the following lemma:

\begin{lemma}
The infinitesimal intersection multiplicity of $\Gamma_s$ at a point $(x,y)$ is equal to
the intersection multiplicity of the curves $\Gamma_s=\{f_s=0\}$ and $\{\frac{\partial f_s}{\partial s}=0\}$
at the same point.
\end{lemma}

This is a direct generalization of Lemma 4.1.

\proof
The Jacobian of $F$ is, by definition, $\det(\pd Fs,\pd Ft)$.
We have $df_s(\pd Ft)=0$, and neither the 1-form $df_s$ nor the
vector field $\pd Ft$ ever vanish.
It follows that the Jacobian of $F$ is $df_s(\pd Fs)$ times a
nowhere vanishing differentiable function.
In particular, the order of vanishing of the Jacobian
coincides with the order of vanishing of the function $df_s(\pd Fs)$ at the same point,
and, by the equality
$$
\pd {f_s}s\circ F+df_s\left(\pd Fs\right)=0,
$$
with the order of vanishing of $\pd {f_s}s\circ F$ at the same point.
Restrict all functions considered to a curve $\Gamma_s$ and express them
it terms of the local parameter $t$.
Then the order of vanishing of the function $\pd {f_s}s$ is, by definition, the
intersection multiplicity of the curves $\Gamma_s$ and $\{\pd {f_s}s=0\}$.
$\Box$

Theorem \ref{Bez} now follows.

The following statement provides a description of families $F$ that consist of
osculating algebraic curves to a given plane curve:

\begin{theorem}
Under the assumptions of Theorem \ref{Bez},
suppose also that there is a smooth plane curve $\gamma$ parameterized by $s$ and
such that $\gamma(s)\in\Gamma_s$ for each $s$, and each curve $\Gamma_s$
has infinitesimal intersection multiplicity $n(d)-1$ at the point $\gamma(s)$.
Then $\Gamma_s$ are osculating algebraic curves of degree $d$ for the curve $\gamma$.
\end{theorem}

This theorem generalizes the well-known algorithm of finding the envelope
of a family of lines: the envelope coincides with the locus of points,
where two infinitesimally close lines intersect.

\proof
Since $\pd{f_s}s=0$ on $\gamma$, the curve $\gamma$ is the envelope of curves $\Gamma_s$
(this follows from the classical description of the envelope).

Then the function $s'\mapsto f_s(\gamma(s'))$ has a multiple zero at point $s'=s$.
By the preparation theorem for differentiable functions \cite{Loj,Hor,Mal},
we have
$$
f_s(\gamma(s'))=(s-s')^m\phi(s,s'),
$$
where $m>1$ is an integer and $\phi$ is a smooth function of two variables such that
$\phi(s,s)\ne 0$ locally near a given value of $s$.
In particular, the curves $\Gamma_s$ approximate the curve $\gamma$ up to
order $m$ for $s$ in the chosen neighborhood.
Reparameterize curves $\Gamma_s$ to make $\Gamma_s(t)$ coincide with $\gamma(s)$
for $t=0$.
Then $\pd{f_s}s(\Gamma_s(t))$ vanishes at point $t=0$ with order $m-1$.
On the other hand, the order of vanishing is $n(d)-1$, hence $m=n(d)$.
$\Box$

\bigskip

{\bf EndNote.}
After this paper was completed, we learned that E.
Ghys has been extensively lecturing on the same subject in the
recent years.
We are pleased to acknowledge his priority in this matter;
some of his unpublished results generalize those given in this paper.
E. Ghys and the authors plan to further develop these ideas
and to publish a more complete and more general account.

\medskip

{\bf Acknowledgments.} We are grateful to Dan Genin for discussions
and for making figures for this paper in Mathematica.
Many thanks to E. Ghys for historical comments, in particular, for attracting out attention to Tait's paper \cite{Tai}.

\end{document}